\documentclass[12pt]{amsart}
\usepackage{amssymb}
\usepackage{graphicx}
\usepackage{amssymb,latexsym,amsthm}

\usepackage{amsmath,amssymb}
\usepackage{amssymb,latexsym,amsthm,mathptmx}
\usepackage{amsmath,amsfonts,amssymb}
\usepackage{graphicx}
\usepackage{color}
\usepackage{mathrsfs}
\usepackage{amsmath,amsfonts,amsthm,amssymb}
\usepackage{newtxtext}
\usepackage{newtxmath}
\usepackage[mathscr]{eucal}
\usepackage{bm}

\numberwithin{equation}{section}

\newtheorem{thm}{Theorem}[section]
\newtheorem{prop}[thm]{Proposition}

\newtheorem{cor}[thm]{Corollary}
\newtheorem{defn}[thm]{Definition}
\newtheorem{exam}[thm]{Example}

\newcommand{\R}{\mathbb{R}}

\newcommand\V{\mathbb V}
\newcommand{\uv}{\mbox{\boldmath $u$}}
\newcommand{\vvv}{\mbox{\boldmath $v$}}

\newcommand\gyr{\operatorname {gyr}}
\newcommand{\e}{\mbox{\boldmath $e$}}
\newcommand{\ab}{\mbox{\boldmath $a$}}
\newcommand{\bb}{\mbox{\boldmath $b$}}
\newcommand{\ccb}{\mbox{\boldmath $c$}}

\newcommand{\x}{\mbox{\boldmath $x$}}
\newcommand{\y}{\mbox{\boldmath $y$}}
\newcommand{\z}{\mbox{\boldmath $z$}}

\newcommand{\PP}{\mbox{\boldmath $P$}}

\begin{document}
\title
{Generalized gyrovector spaces revisited}

\author{Toshikazu \textsc{Abe}
}
\address{Faculty of Engineering, Ibaraki University, Hitachi 316-8511, Japan.}
\email{toshikazu.abe.bin@vc.ibaraki.ac.jp}

\author{Osamu \textsc{Hatori}}
\address{Institute of Science and Technology, Niigata University, Niigata 950-2181, Japan.}
\email{hatori@math.sc.niigata-u.ac.jp}


\begin{abstract}
We revise a proof of a Mazur-Ulam theorem for generalized gyrovector spaces.
\end{abstract}

\maketitle

\section{Introduction}
Ungar introduced the concept of a gyrovector space, which is discussed in detail in \cite{ungyrov,Un}. In contrast, we later proposed the idea of a generalized gyrovector space (GGV) in \cite{abehatori} and proved a Mazur-Ulam theorem for it. However, we have identified some errors in that paper and have revised them in this paper.

We want to clarify that statements of Theorem 13 and Corollary 14 in \cite{abehatori} are correct, but their proofs need some revisions. Specifically, we need to correct an equation on the sixth line of page 398, which should read as $\|\phi(e)\|=0_{\|\phi(G)\|}$ instead of $\|\phi(e)\|=0$. It is important to note that $0_{\|\phi(G)\|}$ may not always represent the real number 0, and we have provided an example in Example \ref{pathological}. Additionally, we would like to mention that all examples of GGVs presented in \cite{abehatori} and \cite{hatori} satisfy $\|\phi(e)\|=0$ (i.e., the real number 0). 

Let $n(\x)=\|\phi(\x)\|$ be the original length on a GGV $G$.
With this length, $G$ is a normed gyrolinear space in the sense of \cite[Definition~2.2]{abe} if and only if $n(\e)=0$.
For every GGV, the map $f$ in Proposition~\ref{lin} gives a normalized length $N=f\circ n$.
With $N$ as its length and the identity map on $[0,\infty)$ as its auxiliary function, $G$ is a normed gyrolinear space.

Equation~\eqref{scin} corrects (10) of \cite{abehatori}. We give revised proofs of Theorem~\ref{GMU} and Corollary~\ref{gl}.

We follow the terminology of \cite{abehatori}. See also \cite{ungyrov,Un}.
\section{Preliminary}
Before recalling the definition of a generalized gyrovector space, it would be convenient to confirm the following. 
\begin{prop}\label{prGGV}
Let $(G,\oplus)$ be a gyrocommutative gyrogroup with the map $\otimes:\R\times G\to G$. 
Let $\e$ denote the identity of $(G,\oplus)$.
Let $\phi$ be an injection from $G$ into a real normed space $(\V,\|\cdot\|)$. 
We assume
\begin{itemize}


\item[{\rm{(GGV1)}}] $1\otimes\ab=\ab$ for any $\ab\in G$; 

\item[{\rm{(GGV2)}}] $(r_1+r_2)\otimes\ab =(r_1\otimes\ab)\oplus (r_2\otimes\ab)$ for any $\ab\in G, \ r_1,r_2\in\R$;

\item[{\rm{(GGV3)}}] $(r_1r_2)\otimes\ab=r_1\otimes(r_2\otimes\ab)$ for any $\ab\in G, \ r_1,r_2\in\R$;

\item[{\rm{(GGVV)}}] $\|\phi(G)\|=\{\pm\|\phi(\ab)\|\in\R:\ab\in G\}$ is a 
real linear space of dimension one with vector addition $\oplus'$ and scalar multiplication $\otimes'$;

\item[{\rm{(GGV7)}}] $\|\phi(r\otimes\ab)\|=|r|\otimes'\|\phi(\ab)\|$ for any $\ab\in G$, $r\in \R$.
\end{itemize}
Then we have the following.
\begin{itemize}
\item[{\rm{(i)}}]
$\|\phi(\e)\|=0_{\|\phi(G)\|}$;
\item[{\rm{(ii)}}]
 $r\otimes \e=\e$ for all $r\in \R$;
\item[{\rm{(iii)}}]
$0\otimes \ab=\e$ for all $\ab\in G$;
\item[{\rm{(iv)}}]
$r\otimes \ab=\e$ for $r\in \R$ and $\ab\in G$ implies that $r=0$ or $\ab=\e$;
\item[{\rm{(v)}}]
$\ominus (\alpha\otimes \ab)=(-\alpha)\otimes \ab$ for all $\alpha\in \R$, $\ab\in G$, in particular, $\ominus \ab=(-1)\otimes \ab$ for all $\ab\in G$;
\item[{\rm{(vi)}}]
$\|\phi(\ab)\|\ne0$ for all $\ab\in G\setminus\{\e\}$, 
$\|\phi(r\otimes\ab)\|\ne 0$ for all $r\in \R\setminus\{0\}$ and $\ab\in G\setminus\{\e\}$.
\end{itemize}
\end{prop}
\begin{proof}
By the condition (GGV2) we have
\[
0\otimes \ab=(0+0)\otimes \ab=0\otimes \ab \oplus 0\otimes \ab
\]
for every $\ab\in G$. 
By the left cancelation law for gyrogroups, we infer that $\e=0\otimes \ab$, proving (iii). In particular,  since $\e=0\otimes \e$ we obtain
\[
\|\phi(\e)\|=\|\phi(0\otimes \e)\|=0\otimes'\|\phi(\e)\|=0_{\|\phi(G)\|}
\]
by (GGV7), proving (i).

For every $r\in \R$ we have by (GGV3) that
\[
r\otimes \e=r\otimes (0\otimes \e)=(r\times 0)\otimes \e=0\otimes \e=\e,
\]
which proves (ii).

Suppose that $r\otimes \ab=\e$ for some $r\in \R$ and $\ab\in G$. If $r\ne 0$, then we have by (GGV1) and (GGV3) that 
\[
\ab=1\otimes \ab=\frac1r\otimes(r\otimes \ab)=\frac1r \otimes\e=\e.
\]
Thus, (iv) holds. 

By (GGV2) we have
\[
(-\alpha)\otimes \ab\oplus \alpha\otimes \ab=(-\alpha+\alpha)\otimes \ab=0\otimes \ab=\e.
\]
The uniqueness of $\ominus (\alpha\otimes \ab)$ ensures that 
$(-\alpha)\otimes \ab=\ominus (\alpha\otimes \ab)$, 
proving (v).

Suppose that $\|\phi(\ab)\|=0$ for some $\ab\in G$. 
Then we obtain
\[
\|\phi(\ominus \ab)\|=\|\phi((-1)\otimes \ab)\|=|-1|\otimes'\|\phi(\ab)\|=\|\phi(\ab)\|=0.
\]
Thus we have $\phi(\ominus \ab)=\phi(\ab)=0_{\V}$ 
as $\V$ is a normed linear space. As $\phi$ is an injection we infer that $\ominus \ab=\ab$, so 
\[ 
\ab=1\otimes \ab=\frac12\otimes(2\otimes \ab)
=\frac12\otimes((1+1)\otimes \ab)
=\frac12\otimes(\ab\oplus \ab) 
=\frac12\otimes(\ominus \ab\oplus \ab)=\frac12\otimes \e=\e.
\]
We conclude that $\|\phi(\ab)\|\ne 0$ if  $\ab \ne \e$. By (iv), we infer that $\|\phi(r\otimes \ab)\|\ne 0$ if $\ab \in G\setminus \{\e\}$ and $r\in \R\setminus \{0\}$. We have proved (vi).

\end{proof}
Note that the proof of Proposition \ref{prGGV} shows that $\|\phi(\ab)\|=0$ only if $\ab=\e$. One might expect $\|\phi(\e)\|=0$ for a GGV, but Example \ref{pathological} illustrates that this is not always the case.

 In this paper we simply write $\alpha\otimes \ab \oplus \beta\otimes\bb$ instead of  $(\alpha\otimes \ab) \oplus (\beta\otimes\bb)$, and $\ominus\alpha\otimes\ab$ instead of $\ominus(\alpha\otimes\ab)$ for any $\alpha,\beta\in\R$, $\ab,\bb\in G$. 
\section{Definition of a generalized gyrovector space revisited}
We recall the definition of a generalized gyrovector space (cf. \cite{abehatori}).
\begin{defn}[A generalized gyrovector space]\label{GGV}
Let $(G,\oplus)$ be a gyrocommutative gyrogroup with the map $\otimes:\R\times G\to G$. 
Let $\e$ denote the identity of $(G,\oplus)$. 
Let $\phi$ be an injection from $G$ into a real normed space $(\V,\|\cdot\|)$. 
We say that $(G,\oplus,\otimes,\phi)$ (or $(G,\oplus,\otimes)$ just for a simple notation) is a generalized gyrovector space or a GGV in short if  the following conditions (GGV0) to (GGV8) are fulfilled:
\begin{itemize}

\item[{\rm{(GGV0)}}] $\|\phi(\gyr[\uv,\vvv]\ab)\|=\|\phi(\ab)\|$ for any $\uv,\vvv,\ab\in G$;

\item[{\rm{(GGV1)}}] $1\otimes\ab=\ab$ for any $\ab\in G$; 

\item[{\rm{(GGV2)}}] $(r_1+r_2)\otimes\ab =(r_1\otimes\ab)\oplus (r_2\otimes\ab)$ for any $\ab\in G, \ r_1,r_2\in\R$;

\item[{\rm{(GGV3)}}] $(r_1r_2)\otimes\ab=r_1\otimes(r_2\otimes\ab)$ for any $\ab\in G, \ r_1,r_2\in\R$;

\item[{\rm{(GGV4)}}] $\phi(|r|\otimes\ab) / \|\phi( r\otimes\ab)\|=\phi(\ab) / \|\phi(\ab)\|$ for any $\ab\in G\setminus\{\e\}, \ r\in\R\setminus\{0\}$, where $\e$ denotes the identity element of the gyrogroup $(G,\oplus)$;

\item[{\rm{(GGV5)}}] $\gyr[\uv,\vvv](r\otimes\ab)=r\otimes \gyr[\uv,\vvv]\ab$ for any $\uv,\vvv,\ab\in G, \ r\in\R$;

\item[{\rm{(GGV6)}}] $\gyr[r_1\otimes\vvv,r_2\otimes\vvv]=id_G$  for any $\vvv\in G, \ r_1,r_2\in\R$, where $id_G$ is the identity map on $G$;

\item[{\rm{(GGVV)}}] $\|\phi(G)\|=\{\pm\|\phi(\ab)\|\in\R:\ab\in G\}$ is a 
real linear space of dimension one with vector addition $\oplus'$ and scalar multiplication $\otimes'$; 

\item[{\rm{(GGV7)}}] $\|\phi(r\otimes\ab)\|=|r|\otimes'\|\phi(\ab)\|$ for any $\ab\in G$, $r\in \R$;

\item[{\rm{(GGV8)}}] $\|\phi(\ab\oplus\bb)\|\leq\|\phi(\ab)\|\oplus'\|\phi(\bb)\|$ for any $\ab,\bb \in G$.
\end{itemize}
\end{defn}
The quotients in (GGV4) are well defined since Proposition~\ref{prGGV}(vi) gives
$\|\phi(\ab)\|\ne 0$ and $\|\phi(r\otimes \ab)\|\ne 0$ for every $\ab\in G\setminus \{\e\}$ and $r\in \R\setminus\{0\}$.

\begin{prop}\label{prGGV2}
Suppose that $(G,\oplus,\otimes,\phi)$ is a GGV. Then 
$\|\phi(\ab)\|=0_{\|\phi(G)\|}$ if and only if $\ab=\e$, where $\e$ is the unit of $G$. We have 
$r\otimes' \|\phi(\ab)\|=s\otimes' \|\phi(\ab)\|$ if and only if $r=s$ or $\ab=\e$. 
\end{prop}
\begin{proof}
Suppose that $\ab=\e$. Then by (i) of Proposition \ref{prGGV} we have $\|\phi(\ab)\|=0_{\|\phi(G)\|}$.

Suppose conversely that $\|\phi(\ab)\|=0_{\|\phi(G)\|}$. We prove $\ab=\e$. 
Suppose not; $\ab\ne \e$. Then we have by (GGV4) that 
\[
\frac{\phi(2\otimes \ab)}{\|\phi(2\otimes \ab)\|}=\frac{\phi(\ab)}{\|\phi(\ab)\|}.
\]
On the other hand 
\[
\|\phi(2\otimes \ab)\|=2\otimes'\|\phi(\ab)\|=2\otimes'0_{\|\phi(G)\|}=0_{\|\phi(G)\|}=\|\phi(\ab)\|
\]
ensure that 
\[
\phi(2\otimes \ab)=\phi(\ab).
\]
As $\phi$ is injective, we infer that $\ab\oplus \ab=2\otimes \ab=\ab$, so that $\ab=\e$ by the left cancelation law, which is a contradiction since we assume that $\ab\ne \e$. Thus, we have that 
$\ab=\e$.

We prove the second statement. 
Suppose that $r=s$. Then $r\otimes' \|\phi(\ab)\|=s\otimes' \|\phi(\ab)\|$ follows. Suppose that $\ab=\e$. Then $\|\phi(\ab)\|=0_{\|\phi(G)\|}$, so that $r\otimes' \|\phi(\ab)\|=s\otimes' \|\phi(\ab)\|$ follows. 

Suppose that $r\otimes' \|\phi(\ab)\|=s\otimes' \|\phi(\ab)\|$. Then $r=s$ 
or $\|\phi(\ab)\|=0_{\|\phi(G)\|}$. It follows that $r=s$ or $\ab=\e$ by the first part.
\end{proof}

\begin{exam}\label{pathological}
Let $G=(-\infty, -1)\cup [1, \infty)$. Put a function $\Phi\colon \R\to G$ defined by
\begin{equation*}
\Phi(x) =
\begin{cases}
\exp x,\qquad & x\ge 0 \\
-\exp(-x), \qquad & x<0.
\end{cases}
\end{equation*}
Then $\Phi$ is a bijection. Define
\[
\ab\oplus \bb=\Phi(\Phi^{-1}(\ab)+\Phi^{-1}(\bb)),\qquad \ab, \bb\in G.
\]
Then the additive group structure of $\R$ is transplanted to $G$ through $\Phi$. Hence $(G,\oplus)$ is a commutative group. 
Thus $(G,\oplus)$ is a gyrocommutative gyrogroup such that the gyration $\gyr[\ab, \bb]$ is the identity for all $\ab, \bb\in G$. 
The map $\otimes \colon \R\times G\to G$ is defined by 
\[
r\otimes \ab=\Phi(r\Phi^{-1}(\ab)), \qquad r\in \R, \ab\in G.
\]
It is easy to see that (GGV1), (GGV2), (GGV3), (GGV5) and (GGV6) hold. 
Let $(\V,\|\cdot\|)=(\R,|\cdot|)$. 
We define the map $\phi\colon G\to \V$ by
\[
\phi(\ab)=\ab, \qquad \ab\in G.
\]
Then $\|\phi(\ab)\|=|\ab|$ for every $\ab\in G$. It is clear that (GGV0) holds since $\gyr[\cdot,\cdot]$ is the identity. 
By a simple calculation, we have $\|\phi(r\otimes \ab)\|=|\ab|^{|r|}$ 
for all $\ab\in G$ and $r\in \R$. 
For $r\ne 0$, we also obtain $\phi(|r|\otimes \ab)=\ab^{|r|}$ if $\ab\ge 1$, $\phi(|r|\otimes \ab)=-|\ab|^{|r|}$ if $\ab<-1$. Hence (GGV4) holds. 
Let 
\[
\|\phi(G)\|_+=\{\|\phi(\ab)\|\colon \ab\in G\},
\]
\[
\|\phi(G)\|_-=\{-\|\phi(\ab)\|\colon \ab\in G\}. 
\]
We infer $\|\phi(G)\|_+=[1,\infty)$, $\|\phi(G)\|_-=(-\infty, -1]$ and $\|\phi(G)\|_+\cup \|\phi(G)\|_-=\|\phi(G)\|$.  
Suppose that $S\colon (-\infty,0)\to (-\infty, -1]$ is a bijection. Let $T\colon \R\to \|\phi(G)\|$ be a bijection defined as 
\begin{equation*}
    T(x) =
    \begin{cases}
      \exp x, \quad& x\ge 0 
      \\
      S(x), \quad& x<0.
    \end{cases}
\end{equation*}
Then, $T([0,\infty))=\|\phi(G)\|_+$, $T((-\infty,0))=\|\phi(G)\|_-$.
We define operations $\oplus'$ and $\otimes'$ on $\|\phi(G)\|$ as follows. 
\[
A\oplus' B=T(T^{-1}(A)+T^{-1}(B)),\quad A,B\in \|\phi(G)\|,
\]
\[
r\otimes' A=T(rT^{-1}(A)),\quad r\in \R,\,\, A\in \|\phi(G)\|.
\]
Then $(\|\phi(G)\|,\oplus',\otimes')$ is a real-linear space of dimension one since the algebraic structure is simply transplanted from the usual real linear space $\R$ of dimension one through the map $T$. 
By a simple calculation, we have 
\[
\|\phi(r\otimes \ab)\|=|\ab|^{|r|}=|r|\otimes'\|\phi(\ab)\|, \quad r\in \R, \,\,\ab\in G.
\]
Thus (GGV7) holds. To verify (GGV8), put $u=\Phi^{-1}(\ab)$ and
$v=\Phi^{-1}(\bb)$. Since $|\Phi(t)|=\exp(|t|)$ for every $t\in\R$,
the triangle inequality on $\R$ gives
\[
\begin{aligned}
 \|\phi(\ab\oplus\bb)\|
 &= |\Phi(u+v)|=\exp(|u+v|)\\
 &\le \exp(|u|)\exp(|v|)=|\ab|\,|\bb|
  =\|\phi(\ab)\|\oplus'\|\phi(\bb)\|,
\end{aligned}
\]
where the last equality follows from the definition of $\oplus'$ on $[1,\infty)$.
We conclude that  $(G,\oplus, \otimes, \phi)$ is a GGV, where $\e=1$ and $0_{\|\phi(G)\|}=\|\phi(\e)\|=1$. 
\end{exam}

We recall 
the gyrometric for a GGV which is defined in \cite[Definition 6]{abehatori}.
\begin{defn}[gyrometric]
Let $(G,\oplus,\otimes,\phi)$ be a GGV. Put
\[
\varrho(\ab,\bb)=\|\phi(\ab\ominus\bb)\|,\qquad \ab,\bb\in G.
\]
We call $\varrho(\cdot,\cdot)$ a gyrometric on $(G,\oplus,\otimes,\phi)$.
\end{defn}
For the convenience of the readers, we state and prove some elementary properties of the gyrometric  (as (4)  and 
(8) in \cite{abehatori}).
\begin{prop}\label{conveni}
\[
\varrho(\x\oplus \ab, \x\oplus\bb)=\varrho(\ab,\bb)=\varrho(\ominus\ab,\ominus \bb)=\varrho(\bb,\ab),\qquad \x, \ab,\bb \in G.
\]
\end{prop}
\begin{proof}
By \cite[Theorem 3.13]{Un} we have $(\x\oplus\ab)\ominus(\x\oplus\bb)=\gyr[\x,\ab](\ab\ominus\bb)$, hence
\[
\varrho(\x\oplus \ab,\x\oplus\bb)=\|\phi((\x\oplus\ab)\ominus(\x\oplus\bb))\|=\|\phi(\gyr[\x,\ab](\ab\ominus\bb))\|,
\]
by (GGV0) 
\[
=\|\phi(\ab\ominus\bb)\|.
\]
As $G$ is gyrocommutative, we have
\[
\varrho(\ab,\bb)=\|\phi(\ab\ominus \bb)\|=\|\phi(\gyr[\ab,\ominus\bb](\ominus \bb\oplus \ab))\|,
\]
then by (GGV0), 
\[
=\|\phi(\ominus\bb\oplus\ab)\|=\varrho(\ominus\bb,\ominus\ab).
\]
We also have by (GGV7) and Proposition \ref{prGGV} (v) that
\[
\varrho(\ab,\bb)=\|\phi(\ab\ominus\bb)\|=\|\phi((-1)\otimes(\ab\ominus\bb))\|=\|\phi(\ominus(\ab\ominus\bb))\|,
\]
then by the gyroautomorphic inverse property \cite[Theorem 3.2]{Un}, and the gyrocommutativity, 
and (GGV0)
\[
=\|\phi(\ominus \ab\oplus \bb)\|
=\|\phi(\gyr[\ominus \ab,\bb](\bb\ominus\ab))\|=\|\phi(\bb\ominus \ab)\|=\varrho(\bb,\ab).
\]
\end{proof}

\begin{prop}\label{triangle}
\begin{equation}\label{gyrotriangle}
\varrho(\ab,\bb)\le \varrho(\ab,\ccb)\oplus'\varrho(\ccb,\bb), \qquad \ab,\bb, \ccb \in G.
\end{equation}
\end{prop}
\begin{proof}
By Proposition \ref{prGGV} (v), Proposition \ref{conveni}, (GGV7), and (GGV8), we have
\begin{multline*}
\varrho(\ab,\bb)=\varrho(\ominus \ab,\ominus \bb)=\varrho(\ccb\ominus \ab,\ccb\ominus \bb)
=\|\phi((\ccb\ominus\ab)\ominus(\ccb\ominus \bb))\|
\\
=\|\phi((\ccb\ominus\ab)\oplus(\ominus(\ccb\ominus \bb)))\|
\le \|\phi(\ccb\ominus\ab)\|\oplus' \|\phi(\ominus(\ccb\ominus \bb))\|
\\
=\|\phi(\ccb\ominus \ab)\|\oplus'
(|-1|\otimes'\|\phi(\ccb\ominus\bb)\|)
=\|\phi(\ccb\ominus\ab)\|\oplus' \|\phi(\ccb\ominus \bb)\|
\\
=\varrho(\ccb,\ab)\oplus' \varrho(\ccb,\bb)
=\varrho(\ab,\ccb)\oplus' \varrho(\ccb,\bb).
\end{multline*}
\end{proof}
We call the equation \eqref{gyrotriangle} the gyrotriangle inequality for a GGV.
\begin{defn}[gyrometric preserving map]\label{GMP}
Let $(G_1,\oplus_1,\otimes_1)$ and $(G_2,\oplus_2,\otimes_2)$ be GGVs. 
Let $\varrho_1$ and $\varrho_2$ be gyrometrics of $G_1$ and $G_2$, respectively. 
We say that a map $T:G_1{\to}G_2$ is gyrometric preserving  if the equality
\[\varrho_2(T\ab,T\bb)=\varrho_1(\ab,\bb)\]
holds for every pair $\ab,\bb\in G_1$.
\end{defn}
\section{A Mazur-Ulam theorem for GGVs}
The section presents a theorem stating that a surjective gyrometric preserving map preserves gyromidpoints. Gyromidpoints are defined in Definition \ref{gyromidpoint}. This theorem and its corollary generalize the Mazur-Ulam theorem for real normed spaces.
\begin{defn}[gyrogroup coaddition]\label{coaddition}
    Let $(G,\oplus)$ be a gyrogroup. The gyrogroup coaddition $\boxplus$ is defined by 
    \[
    \ab\boxplus \bb=\ab\oplus \gyr[\ab,\ominus \bb]\bb,\quad \ab,\bb\in G.
    \]
\end{defn}
It is known that the gyrogroup coaddition $\boxplus$ is commutative if and only if the gyrogroup is gyrocommutative \cite[Theorem 3.4]{Un}.
\begin{defn}[gyromidpoint]\label{gyromidpoint}
Let $(G,\oplus,\otimes)$ be a GGV.
The gyromidpoint $\PP(\ab,\bb)$ of $\ab,\bb\in (G,\oplus,\otimes)$ is defined as $\PP(\ab,\bb)=\frac{1}{2}\otimes(\ab\boxplus\bb)$, where $\boxplus$ is the gyrogroup coaddition of the gyrogroup $(G,\oplus)$. 
\end{defn}
Note that $\PP(\ab,\bb)=\PP(\bb,\ab)$ as  $\boxplus$ is commutative.
The argument in the proof of \cite[Theorem~6.34]{Un} gives
\begin{equation}\label{gme}
\PP(\ab,\bb) = \ab\oplus\frac{1}{2}\otimes(\ominus\ab\oplus\bb)
\end{equation}
(see also \cite[Definition~6.32]{Un}).
In particular, $\PP(\ab,\bb)=\frac{1}{2}(\ab + \bb)$ if the GGV $(G,\oplus,\otimes)$ is indeed a real normed space $(G,+,\cdot)$. 
The equation (9) of  Proposition 15 in \cite{abehatori} is
\begin{equation}\label{gyromid}
\varrho(\ab,\PP(\ab,\bb))=\varrho(\bb, \PP(\ab,\bb))=\frac12 \otimes' \varrho(\ab,\bb).
\end{equation}
The following appears as Proposition 17 in \cite{abehatori}. We can similarly prove this as the proof of \cite[Theorem 2.28]{Un}. 
\begin{prop}\label{iso}
Let $(G_1,\oplus_1)$ and $(G_2,\oplus_2)$ be gyrogroups. 
Suppose that $T: G_1{\to}G_2$ is a bijection. 
Then, {\rm{(I1)}} and {\rm{(I2)}} are equivalent to each other. 
\begin{description}
\item[{\rm{(I1)}}] $T(\ab\oplus_{1}\bb)=T(\ab)\oplus_2T(\bb)$ for any $\ab,\bb \in G_1$,

\item[{\rm{(I2)}}] $T(\ab\boxplus_{1}\bb)=T(\ab)\boxplus_2T(\bb)$ for any $\ab,\bb \in G_1$.
\end{description}
\end{prop}
\begin{thm}[Theorem 13 in 
 \cite{abehatori}]\label{GMU}
Let $(G_1,\oplus_1,\otimes_1)$ and $(G_2,\oplus_2,\otimes_2)$ be GGVs. 
Let $\varrho_1$ and $\varrho_2$ be gyrometrics of $G_1$ and $G_2$, respectively.
Suppose that a surjection $T\colon G_1{\to}G_2$ satisfies 
\[\varrho_2(T\ab,T\bb)=\varrho_1(\ab,\bb)\]
 for any pair $\ab,\bb\in G_1$. 
Then $T$ preserves the gyromidpoints:
\[\PP(T\ab,T\bb) =T(\PP(\ab,\bb))\]
for any pair $\ab,\bb\in G_1$.
\end{thm}
The next corollary states that a surjective map that preserves the gyrometric structure also preserves the algebraic structure modulo the left gyrotranslations.
\begin{cor}[Corollary 14 in \cite{abehatori}]\label{gl}
Let $(G_1,\oplus_1,\otimes_1)$ and $(G_2,\oplus_2,\otimes_2)$ be GGVs. 
Let $\varrho_1$ and $\varrho_2$ be gyrometrics of $G_1$ and $G_2$, respectively. 
Suppose that a surjection $T:G_1{\to}G_2$ satisfies 
\[\varrho_2(T\ab,T\bb)=\varrho_1(\ab,\bb)\]
 for any pair $\ab,\bb\in G_1$. 
Then $T$ is of the form $T=T(\e_1)\oplus_2 T_0$, 
where $T_0$ is an isometrical isomorphism in the sense that it is a bijection such that the equalities
\begin{eqnarray}
&&T_0(\ab\oplus_1 \bb) = T_0(\ab) \oplus_2 T_0(\bb); \label{c1} \\
&&T_0(\alpha \otimes_1 \ab) = \alpha \otimes_2 T_0(\ab); \label{c2} \label{6}\\
&&\varrho_2(T_0(\ab),T_0(\bb))=\varrho_1(\ab,\bb) \label{c3}
\end{eqnarray}
for all  $\ab,\bb \in G_1$ and $\alpha \in \R$ hold.
\end{cor}

\section{Preparations of the proof}\label{lemma}

Let $(G,\oplus,\otimes,\phi)$ be a GGV. 
Suppose that $\ab\in G$ and $0\leq \alpha \leq \beta$. 
Then 
\begin{eqnarray}\label{sm}
\alpha\otimes'\|\phi(\ab)\| &=& \|\phi(\alpha\otimes\ab)\| \\
                                     &=& \left\|\phi\left(\left( \frac{\beta + \alpha}{2}-\frac{\beta - \alpha}{2}\right)\otimes\ab\right)\right\| \nonumber \\
                                     &=& \left\|\phi\left(\left( \frac{\beta + \alpha}{2}\right)\otimes\ab\oplus\left(-\frac{\beta - \alpha}{2}\right)\otimes\ab\right)\right\| \nonumber \\
                                     &\leq& \left( \frac{\beta + \alpha}{2}\right)\otimes'\|\phi(\ab)\|\oplus'\left(\frac{\beta - \alpha}{2}\right)\otimes'\|\phi(\ab)\| \nonumber \\
                              &=& \left(\frac{\beta + \alpha}{2}+\frac{\beta - \alpha}{2}\right)\otimes'\|\phi(\ab)\| \nonumber \\
                              &=& \beta\otimes'\|\phi(\ab)\| \nonumber
\end{eqnarray}
If $\alpha\ge 0$, then $\alpha\otimes'\|\phi(\ab)\|=\|\phi(\alpha\otimes \ab)\|\ge 0$.
By \eqref{sm} 
and Propositions \ref{prGGV} (vi) and \ref{prGGV2} we have that 
$0\le\alpha<\beta$ implies that 
\[
0\le\|\phi(\alpha\otimes \ab)\|=\alpha\otimes'\|\phi(\ab)\|<\beta\otimes'\|\phi(\ab)\|
\]
for any $\ab\in G\setminus\{\e\}$. Suppose conversely that 
\[
0\le\alpha\otimes'\|\phi(\ab)\|<\beta\otimes'\|\phi(\ab)\|.
\]
for some $\alpha\ge0, \beta\ge 0$ and $\ab\in G\setminus \{\e\}$. 
We prove $\alpha<\beta$. Suppose not: $\alpha\ge \beta$. 
By \eqref{sm} we have 
$\beta\otimes'\|\phi(\ab)\|\le \alpha\otimes'\|\phi(\ab)\|$, which is a contradiction proving 
$\alpha<\beta$.
Therefore, we have 
\begin{equation}\label{scin}
0\le \alpha<\beta \iff 0\le\alpha\otimes'\|\phi(\ab)\|<\beta\otimes'\|\phi(\ab)\|
\end{equation}
for any $\alpha, \beta \ge 0$ and $\ab\in G\setminus\{\e\}$. (Note that equation \eqref{scin} is a corrected version of (10) in \cite{abehatori}, which was partly incorrect. In fact, if  $(G,\oplus,\otimes,\phi)$ is a GGV such that $\|\phi(\e)\|\ne 0$ (cf. Example \ref{pathological}), then $0<\|\phi(\e)\|=0_{\|\phi(G)\|}=0\otimes' \|\phi(\ab)\|$ for $\ab\in G\setminus \{\e\}$ by Proposition \ref{prGGV}. 
Proposition 18 in \cite{abehatori} was proved by an application of the statement (10). We revise it as
Proposition \ref{lin}. 
Note also  that $\|\phi(\e)\|=0$ which appears on line 6 on page 398 of \cite{abehatori} is a typo, it should be read as $\|\phi(\e)\|=0_{\|\phi(G)\|}$.
)
The Proposition 18 mentioned in \cite{abehatori} has been revised and corrected as follows. 
\begin{prop}\label{lin}
Let $(G,\oplus,\otimes,\phi)$ be a GGV. 
Then there exists a bijection $f:\|\phi(G)\|\rightarrow \R$ that satisfies the following conditions;  
\begin{itemize}
\item[{\rm{(F1)}}] $f(A\oplus' B)=f(A)+f(B)$ and $f(r\otimes' A)=rf(A)$ for any $A,B \in \|\phi(G)\|$,  $r\in \R$ 
\item[{\rm{(F2)}}] $0\le A<B$ if and only if $0\le f(A)<f(B)$ for $A,B\in \|\phi(G)\|$ with $A\ge0$, $B\ge 0$.
\end{itemize}
\end{prop}
\begin{proof}
By the condition {\rm{(GGVV)}}, there exists a bijection $f:\|\phi(G)\|\rightarrow \R$ that satisfies the condition {\rm{(F1)}}. 
Needless to say, $-f$ also satisfies the condition {\rm{(F1)}}. 
Hence, we may assume that $f(\|\phi(\x_0)\|)>0$ for some $\x_0\in G$. 

For $0\le A,B\in\|\phi(G)\|$, put $X=\|\phi(\x_0)\|$ and $r_A=f(A)/f(X)$, $r_B=f(B)/f(X)$. 
Then $f(r_A\otimes'X)=r_Af(X)=f(A)$. As $f$ is injective, we infer that $r_A\otimes' X=A$. 
Then we have $r_A\ge 0$. (Suppose not: $r_A<0$. Choose $\ab\in G$ with $\|\phi(\ab)\|=A$,
and put $\bb=(-r_A)\otimes\x_0$ and $z=0_{\|\phi(G)\|}$.
By (GGV7), $\|\phi(\bb)\|=(-r_A)\otimes'X$, so
$A\oplus'\|\phi(\bb)\|=z$. By \eqref{sm} with $\alpha=0$ and
$\beta=1$, and by (GGV8), we have
\[
   z\le \|\phi(\ab\oplus\bb)\|
     \le A\oplus'\|\phi(\bb)\|=z.
\]
Proposition~\ref{prGGV2} gives $\ab\oplus\bb=\e$, hence $\bb=\ominus\ab$.
By Proposition~\ref{prGGV}(v) and (GGV7),
$\|\phi(\bb)\|=\|\phi(\ominus\ab)\|=A$. Thus
$(-r_A)\otimes'X=r_A\otimes'X$. Applying $f$ and using $f(X)>0$
gives $r_A=0$, a contradiction.) In the same way, we have that $r_B\ge 0$.
Then by \eqref{scin} we have
\begin{equation*}
0\le r_A<r_B \iff 0\le A<B
\end{equation*}
since $A=r_A\otimes'X=r_A\otimes'\|\phi(\x_0)\|$ and $B=r_B\otimes'X=r_B\otimes'\|\phi(\x_0)\|$. 
On the other hand, the equations $r_Af(X)=f(A)$ and $r_Bf(X)=f(B)$, with $f(X)>0$, ensure that 
\begin{equation*}
0\le r_A<r_B \iff 0\le f(A)<f(B).
\end{equation*}
Hence we conclude that 
\begin{equation*}
0\le A<B \iff 0\le f(A)<f(B)
\end{equation*}
for every $A,B\in \|\phi(G)\|$ with $A\ge0$ and $B\ge 0$. 
\end{proof}
Note that $f(A)=0$ if and only if $A=0_{\|\phi(G)\|}$ since $f$ is a linear bijection from $\|\phi(G)\|$ onto $\R$. In this case, $A$ need not be $0$ (cf. Example \ref{pathological}). It means that $0<A$ does not always imply $0<f(A)$. 
\begin{prop}\label{clin}
Let $(G,\oplus,\otimes,\phi)$ be a GGV. Suppose that $A,B,A',B'\in \|\phi(G)\|$. 
Then $0\le A<B$ and $0\le A'<B'$ imply that $0\le A\oplus' A'<B\oplus' B'$.
\end{prop}
\begin{proof}
Suppose that $0\le A<B$ and $0\le A'<B'$.
Let $D=\{\|\phi(\x)\|:\x\in G\}$. With $X=\|\phi(\x_0)\|$ as in
the proof of Proposition~\ref{lin}, every element of $D$ has the form
$r\otimes' X$ with $r\ge 0$. Conversely, (GGV7) gives
$r\otimes' X=\|\phi(r\otimes\x_0)\|\in D$ for every $r\ge 0$.
Since $(r\otimes' X)\oplus'(s\otimes' X)=(r+s)\otimes' X$,
$D$ is closed under $\oplus'$. In particular, $A\oplus' A'$ and
$B\oplus' B'$ belong to $D$ and are nonnegative.
By Proposition \ref{lin} {\rm{(F2)}} we have $0\le f(A)<f(B)$  and $0\le f(A')<f(B')$. 
By Proposition \ref{lin} {\rm{(F1)}} we have
\[
0\le f(A\oplus' A')=f(A)+f(A')<f(B)+f(B')=f(B\oplus' B').
\]
Then  Proposition \ref{lin} {\rm{(F2)}} ensures $0\le A \oplus' A'<B\oplus' B'$.

\end{proof}
The gyrometric $\varrho$ need not be a metric on $G$.
For example, in Example~\ref{pathological} we have $\varrho(\ab,\ab)=1$ for every $\ab\in G$.
The map $f$ in Proposition~\ref{lin} gives a metric $f\circ\varrho$ on $G$.
\section{Proof of Theorem \ref{GMU}}\label{proof}
We use Nica's defect argument \cite{nica} to prove Theorem~\ref{GMU}.
We recall Proposition 16 in \cite{abehatori}.
\begin{prop}[Proposition 16 in \cite{abehatori}]\label{ahp16}
    Let $(G,\oplus,\otimes)$ be a GGV. Fix $\ab\in G$ and define $\varphi\colon G\to G$ by $\varphi(\x)=2\otimes \ab\ominus \x$ for $\x\in G$. Then $\varphi$ is a bijection that satisfies the following. 
    \begin{itemize}
        \item[(p1)] $\varphi\circ\varphi$ is the identity map on $G$;
        \item[(p2)] $\varrho (\varphi(\x),\varphi(\y))=\varrho(\x,\y)$, \quad $\x,\y\in G$;
        \item[(p3)] $\varphi(\x)=\x$ if and only if $\x=\ab$;
        \item[(p4)] $\varphi(\x)=\y$ and $\varphi(\y)=\x$ if $\ab=\PP(\x,\y)$;
        \item[(p5)] $\varrho(\varphi(\x),\x)=2\otimes'\varrho(\ab, \x)$
    \end{itemize}
    for all $\x,\y\in G$. 
\end{prop}
\begin{proof}[Proof of Theorem \ref{GMU}]
We first show that $T$ is injective. If $T(\x)=T(\y)$, then
\[
 \varrho_1(\x,\y)=\varrho_2(T(\x),T(\y))
 =\varrho_2(T(\x),T(\x))=\varrho_1(\x,\x)=0_{\|\phi_1(G_1)\|_1}.
\]
By Proposition~\ref{prGGV2}, $\x\ominus_1\y=\e_1$, so $\x=\y$.
Thus $T$ is bijective, since it is surjective by hypothesis.

First, consider the case where $G_1$ and $G_2$ are the same GGV. For the first part of the proof, we denote $(G,\oplus,\otimes)$ instead of $(G_i,\oplus_i,\otimes_i)$ to simplify the notations. 
Let $f$ be the corresponding map for $G$, which appears in Proposition \ref{lin}. 
Let $\x_1,\x_2\in G$ be arbitrary. Put $p=\PP(\x_1,\x_2)$, $p'=\PP(T(\x_1),T(\x_2))$, and the defect: $d=f(\varrho(T(p),p'))$. We prove $d=0$. 
Let $\varphi_1\colon G\to G$ be given by $\varphi_1(\x)=2\otimes p\ominus \x$ for $\x\in G$ and $\varphi_2\colon G\to G$ by $\varphi_2(\z)=2\otimes p'\ominus \z$ for $\z\in G$. 
Put $S=\varphi_1\circ T^{-1}\circ \varphi_2\circ T$.
Then $S\colon G\to G$ is a surjective gyrometric preserving map such that $S(\x_1)=\x_1$ and $S(\x_2)=\x_2$ by Proposition \ref{ahp16} (p2) and (p4).  
We prove 
\begin{equation}\label{translationlemma}
f(\varrho(S^{2^n}(p), p))=2^{n+1}d
\end{equation}
for every nonnegative integer $n$. 
First, we have by Proposition \ref{ahp16} (p1) and (p2) that 
\begin{equation}\label{yoi}
\varrho(\varphi_1\circ S(\x),\y)=\varrho(T^{-1}\circ\varphi_2\circ T(\x), \y)
=
\varrho(\x, T^{-1}\circ\varphi_2\circ T(\y))=\varrho(\varphi_1(\x),S(\y)),\quad \x,\y\in G.
\end{equation}
As $\varphi_1(p)=p$ we have 
\[
\varrho(S(p), p)=\varrho(\varphi_1\circ S(p), p)
=
\varrho(\varphi_2\circ T(p), T(p))=2\otimes'\varrho(T(p),p').
\]
It follows that $f(\varrho(S(p), p))=2d$, which is \eqref{translationlemma} for $n=0$. 
Suppose that \eqref{translationlemma} holds for $n$. 
Since $2\otimes'\varrho(\x,p)=\varrho(\varphi_1(\x), \x)$ for any $\x\in G_1$ by Proposition \ref{ahp16} (p5), we have
\[
2^{n+2}d=f(2\otimes'\varrho(S^{2^n}(p),p))=f(\varrho(\varphi_1\circ S^{2^n}(p), S^{2^n}(p))),
\]
then applying \eqref{yoi} $2^n$ times 
\[
=f(\varrho(\varphi_1(p), S^{2^{n+1}}(p)))=
f(\varrho(p, S^{2^{n+1}}(p))).
\]
as $\varphi_1(p)=p$. Then, by induction, we obtain \eqref{translationlemma} for every nonnegative integer $n$. 
Since $S(\x_1)=\x_1$ we have by the gyrometric triangle inequality that
\[
\varrho(S^{2^n}(p),p)\le \varrho(S^{2^n}(p), S^{2^n}(\x_1))\oplus'\varrho(S^{2^n}(\x_1), p)
=
\varrho(p,\x_1)\oplus'\varrho(\x_1,p)
\]
for every nonnegative integer $n$. 
Thus 
\[
2^{n+1}d=f(\varrho(S^{2^n}(p),p))\le 2f(\varrho(\x_1,p))
\]
for every nonnegative integer $n$. It follows that $f(\varrho(T(p),p'))=d=0=f(0_{\|\phi(G)\|})$. As $f$ is a bijection, we have 
 \[
 \|\phi(T(p)\ominus p')\|=\varrho(T(p), p')=\|\phi(\e)\|=0_{\|\phi(G)\|}.
 \]
Then by Proposition \ref{prGGV2} we see that 
$T(p)\ominus p'=\e$. Hence $T(p)=p'$. 

We consider the general case. 
Suppose that $T\colon G_1{\to}G_2$ satisfies 
\[
\varrho_2(T\ab,T\bb)=\varrho_1(\ab,\bb)
\]
 for any pair $\ab,\bb\in G_1$. 
 Let $\x_1,\x_2\in G_1$ be arbitrary. 
 Put $S=\varphi_1\circ T^{-1}\circ\varphi_2\circ T$,  
 where $\varphi_1\colon G_1\to G_1$ is given by $\varphi_1(\x)=2\otimes_1 p\ominus_1 \x$ for $\x\in G_1$ and $\varphi_2\colon G_2\to G_2$ is given by $\varphi_2(\z)=2\otimes_2p'\ominus_2 \z$ for $\z\in G_2$, where 
 $p=\PP(\x_1,\x_2)$, $p'=\PP(T(\x_1),T(\x_2))$. Then $S\colon G_1\to G_1$ is a surjective gyrometric preserving map from $G_1$ onto itself. By the first part, we obtain that $S(\PP(\x_1,\x_2))=\PP(S(\x_1),S(\x_2))$ for every $\x_1,\x_2\in G_1$. 
 As $S(\x_1)=\x_1$ and $S(\x_2)=\x_2$, we have $S(p)=p$. 
 By this equation, we have $\varphi_2\circ T(p)=T(p)$ since $\varphi_1(p)=p$. As $\varphi_2(\y)=\y$ if and only if $y=p'$, we obtain $T(p)=p'$.

\end{proof}
\section{Proof of Corollary \ref{gl}}\label{proofofCOr}
We prove that a bijective gyrometric preserving map essentially preserves the algebraic structure.

\begin{proof}[Proof of Corollary \ref{gl}]
As shown at the beginning of the proof of Theorem~\ref{GMU}, $T$ is bijective.

If $G_1$ is a singleton, then $G_2$ is also a singleton. In this case, every statement is clear. In the following, we assume $G_1$ is not a singleton. Then $G_2$ is not a singleton as well since $T$ is a bijection.
Let $T_0=\ominus_2 T(\e_1)\oplus_2 T$. Indeed, $T_0:G_1\rightarrow G_2$ is a bijection and $T=T(\e_1)\oplus_2T_0$ by the left
cancelation law, and $T_0(\e_1)=\e_2$. 
By Proposition \ref{conveni}, $T_0$ is a gyrometric preserving map because so is $T$.  
In fact, 
Proposition \ref{conveni} observes that 
\begin{multline*}
\varrho_2(T_0(\ab),T_0(\bb))=\varrho_2(T(\e_1)\oplus_2 T_0(\ab),T(\e_1)\oplus_2 T_0(\bb))
\\
=
\varrho_2(T(\ab),T(\bb))=\varrho_1(\ab,\bb), \quad \ab,\bb\in G_1.
\end{multline*}
Applying Theorem \ref{GMU} to $T_0$, we have 
\begin{equation}\label{1.11}
T_0(\frac{1}{2}\otimes_1 (\ab\boxplus_1 \bb))=\frac{1}{2}\otimes_2 (T_0(\ab)\boxplus_2 T_0(\bb))
\end{equation}
for any $\ab,\bb \ \in G_1$. Since $T_0(\e_1)=\e_2$, we have
\begin{equation}
\begin{split}
T_0(\frac{1}{2} \otimes_1 \x) &= T_0(\frac{1}{2}\otimes_1 (\x\boxplus_1 \e_1)) \\
    &=\frac{1}{2}\otimes_2 (T_0(\x)\boxplus_2 T_0(\e_1)) \\
    &= \frac{1}{2} \otimes_2 T_0(\x) \label{2ot}
\end{split}
\end{equation}
for any $\x \in G_1$. It follows by equations \eqref{1.11}  and \eqref{2ot} that 
\begin{equation}
T_0(\ab\boxplus_1 \bb) = T_0(\ab)\boxplus_2 T_0(\bb) \label{box1}
\end{equation}
for any $\ab,\bb \ \in G_1$. 
Since $T_0$ is bijective,  we have the equation (\ref{c1}) by Proposition \ref{iso}.

Next we show that $T_0$ is homogeneous: $T_0(\alpha \otimes_1 \ab) = \alpha \otimes_2 T_0(\ab)$ for any $\ab \in G_1$ and $\alpha \in \R$. First we consider  $\alpha=\frac{m}{2^n}$ for some integers $m,n$.
For any $\ab \in G_1$ and for any positive  integer $m$, $T_0(m\otimes_1 \ab)=m\otimes_2 T_0(\ab)$ is satisfied by the equation (\ref{c1}). 
We have $\e_2=T_0(\e_1)=T_0(\ab\ominus_1\ab)=T_0(\ab)\oplus_2T_0(\ominus_1\ab)$ for any $\ab\in G_1$. Hence $T_0(\ominus_1\ab)=\ominus_2T_0(\ab)$. It follows by a routine argument that $T_0(m\otimes_1\ab)=m\otimes_2T_0(\ab)$ for every integer $m$ and $\ab\in G_1$. 
By the equation (\ref{2ot}), we have 
\begin{equation}\label{aew}
T_0(\frac{m}{2^n}\otimes_1 \ab) = \frac{m}{2^n}\otimes_2 T_0(\ab)
\end{equation}
for any integers $m,n$. 
Next, we consider a general $\alpha \in \R$. 
To prove the homogeneity of $T_0$ for a general $\alpha$, we will show that 
\begin{equation}\label{base}
\varrho_2(T_0(\alpha \otimes_1 \ab), \alpha \otimes_2 T_0(\ab)) <\|\phi_2(T_0(\bb))\|_2
\end{equation}
for every $\bb\in T_0^{-1}(G_2\setminus\{\e_2\})$. As 
\[
\varrho_2(T_0(\alpha \otimes_1 \ab), \alpha \otimes_2 T_0(\ab))=
\|\phi_2(T_0(\alpha \otimes_1 \ab)\ominus_2\alpha \otimes_2 T_0(\ab))\|_2,
\]
\eqref{base} implies
\[
T_0(\alpha \otimes_1 \ab)\ominus_2\alpha \otimes_2 T_0(\ab)\ne T_0(\bb)
\]
for every $\bb\in T_0^{-1}(G_2\setminus\{\e_2\})$. As $T_0$ is surjective we will arrive at 
\[
T_0(\alpha \otimes_1 \ab)\ominus_2\alpha \otimes_2 T_0(\ab)=\e_2,
\]
so that 
\[
T_0(\alpha \otimes_1 \ab)=\alpha \otimes_2T_0(\ab),
\]
which is the desired equality. We now prove \eqref{base}. 
Let $(r_k)$ be a sequence of dyadic rationals converging to $\alpha$. 
As $T_0(r_k\otimes_1\ab)=r_k\otimes_2T_0(\ab)$ by \eqref{aew} we have
by the gyrotriangle inequality Proposition \ref{triangle} that
\begin{equation}\label{119}
\varrho_2(T_0(\alpha \otimes_1 \ab), \alpha \otimes_2 T_0(\ab)) 
   \leq  \varrho_2(T_0(r_k\otimes_1 \ab), T_0(\alpha \otimes_1 \ab)) \oplus'_2 \varrho_2(r_k \otimes_2 T_0(\ab), \alpha \otimes_2 T_0(\ab)).
\end{equation}
We estimate the terms on the right-hand side of the inequality. 
By Proposition \ref{lin}, there is a bijection $f_i:\| \phi_i(G_i)\|_i\to \R$ which satisfies the conditions {\rm{(F1)}} and {\rm{(F2)}} $(i=1,2)$. 
Let $\bb \in T_0^{-1}(G_2\setminus \{\e_2\})$ be arbitrary. Note that $T_0^{-1}(G_2\setminus \{\e_2\})\ne \emptyset$ since $T_0$ is a surjection and $G_2$ is not a singleton. 
Then  $\bb\ne\e_1$ since $T_0(\e_1)=\e_2$. Hence $\|\phi_1(\bb)\|_1>0$  by Proposition \ref{prGGV} (vi). Then $f_1(\|\phi_1(\bb)\|_1)\ge0$ by Proposition \ref{lin}. In fact, $f_1(\|\phi_1(\bb)\|_1)>0$. (Suppose that  $f_1(\|\phi_1(\bb)\|_1)=0$. Then $\|\phi_1(\bb)\|_1=0_{\|\phi_1(G_1)\|_1}$ since $f_1$ is a linear bijection from $\|\phi_1(G_1)\|_1$ onto $\R$. By Proposition \ref{prGGV2} we have $\bb=\e_1$, which is a contradiction.) As  $T_0(\bb)\ne \e_2$, we also see that $f_2(\|\phi_2(T_0(\bb))\|_2)>0$. 
There is 
a positive integer  $k_0$  such that $k>k_0$ implies that
\[2|r_k-\alpha|f_1(\|\phi_1(\ab)\|_1)<f_1(\|\phi_1(\bb)\|_1)\]
and
\[2|r_k-\alpha|f_2(\|\phi_2(T_0(\ab))\|_2)<f_2(\|\phi_2(T_0(\bb))\|_2).\]  
As $f_1$ and $f_2$ are linear we infer by (GGV7) that
\begin{equation}\label{eeee}
0\le f_1(\|\phi_1((r_k-\alpha)\otimes_1\ab)\|_1)<f_1(\frac12\otimes_1'\|\phi_1(\bb)\|_1)
\end{equation}
and
\begin{equation}\label{eeeee}
0\le f_2(\|\phi_2((r_k-\alpha)\otimes_2T_0(\ab))\|_2)<f_2(\frac12\otimes_2'\|\phi_2(T_0(\bb))\|_2).
\end{equation}
On the other hand, we have
\begin{multline}\label{alll}
\frac12\otimes'_1\|\phi_1(\bb)\|_1=\|\phi_1(\frac12\otimes_1 \bb\ominus_1 \e_1)\|_1=\varrho_1(\frac12\otimes_1\bb,\e_1)
\\
=\varrho_2(T_0(\frac12\otimes_1\bb),T_0(\e_1))
=\varrho_2(\frac12\otimes_2T_0(\bb),\e_2)
\\
=\|\phi_2(\frac12\otimes_2 T_0(\bb))\|_2=\frac12\otimes'_2\|\phi_2(T_0(\bb))\|_2.
\end{multline}
Applying Proposition \ref{lin} (F2) for the inequalities \eqref{eeee} and \eqref{eeeee} we have by \eqref{alll} that
\begin{equation}\label{117}
\|\phi_1((r_k-\alpha)\otimes_1\ab)\|_1<\frac{1}{2}\otimes'_1\|\phi_1(\bb)\|_1=\frac{1}{2}\otimes'_2\|\phi_2(T_0(\bb))\|_2,
\end{equation}
\begin{equation}\label{118}
\|\phi_2((r_k-\alpha)\otimes_2T_0(\ab))\|_2<\frac{1}{2}\otimes'_2\|\phi_2(T_0(\bb))\|_2.
\end{equation}
Applying \eqref{117} and \eqref{118}  we compute
\begin{multline}\label{120}
\varrho_2(T_0(r_k\otimes_1 \ab),T_0(\alpha\otimes_1\ab))=\varrho_1(r_k\otimes_1\ab,\alpha\otimes_1\ab)
\\
=
\|\phi_1(r_k\otimes_1\ab\ominus_1\alpha\otimes_1\ab)\|_1=
\|\phi_1((r_k-\alpha)\otimes_1\ab)\|_1<\frac{1}{2}\otimes'_2\|\phi_2(T_0(\bb))\|_2.
\end{multline}
\begin{multline}\label{121}
\varrho_2(r_k \otimes_2 T_0(\ab), \alpha \otimes_2 T_0(\ab))
=
\|\phi_2(r_k \otimes_2T_0(\ab)\ominus_2\alpha\otimes_2T_0(\ab))\|_2
\\
= \|\phi_2((r_k-\alpha)\otimes_2T_0(\ab))\|_2<\frac{1}{2}\otimes'_2\|\phi_2(T_0(\bb))\|_2
\end{multline}
Applying Proposition \ref{clin} for \eqref{120} and \eqref{121} we have
\begin{multline}\label{124}
 \varrho_2(T_0(r_k\otimes_1 \ab),T_0(\alpha\otimes_1\ab))\oplus_2'\varrho_2(r_k \otimes_2 T_0(\ab), \alpha \otimes_2 T_0(\ab))
\\
<\frac{1}{2}\otimes'_2\|\phi_2(T_0(\bb))\|_2\oplus_2'\frac{1}{2}\otimes'_2\|\phi_2(T_0(\bb))\|_2
=
\|\phi_2(T_0(\bb))\|_2.
\end{multline}
Combining  \eqref{119} and \eqref{124}, we obtain \eqref{base}. 
Therefore we have $T_0(\alpha \otimes_1 \ab)=\alpha \otimes_2T_0(\ab)$ for every $\ab\in G_1$ and $\alpha\in \R$.

\end{proof}

\subsection*{Acknowledgments}
The authors express their sincerest appreciation to the referee for their valuable comments and advice, which have substantially improved the presentation of this paper.
The first author was supported by 
JSPS KAKENHI Grant Numbers JP22K03332.
The second author was supported by 
JSPS KAKENHI Grant Numbers JP19K03536. 

\begin{thebibliography}{99}
\bibitem{abe}
T.~Abe,
{\em 
Normed gyrolinear spaces: a generalization of normed spaces based on gyrocommutative gyrogroups},
Mathematics Interdisciplinary Research 
{\bf 1} (2016), 143--172

\bibitem{abehatori}
T.~Abe and O.~Hatori,
{\em 
{Generalized Gyrovector Spaces and a Mazur-Ulam theorem}},
Publ. Math. Debrecen
{\bf 87} (2015), 393--413;
\emph{
A Note on the Proof of Theorem 13 in the Paper "Generalized Gyrovector Spaces and a Mazur-Ulam theorem"}, arXiv:1511.05187v1

\bibitem{hatori}
O.~Hatori,
{\em
Examples and applications of generalized gyrovector spaces},
Result. Math. {\bf 71}(2017), 295--317

\bibitem{nica}
B.~Nica,
\emph{
The Mazur-Ulam theorem},
Expo. Math. {\bf 30} (2012), 397--398

\bibitem{ungyrov}
A.~A.~Ungar,
\emph{
Analytic Hyperbolic Geometry: Mathematical Foundations and Applications},
World Scientific,
\textbf(2005)

\bibitem{Un}
A. A. Ungar, 
\emph{Analytic Hyperbolic Geometry and Albert Einstein's Special Theory of Relativity},
World Scientific, 
\textbf (2008)






\end{thebibliography}
\end{document}